\documentclass[11pt]{article}
\usepackage[latin1]{inputenc}
\usepackage[activeacute,spanish]{babel}
\usepackage{graphicx}
\usepackage{multicol}
\usepackage{amsmath,amssymb,epsfig}
\usepackage{float}
\usepackage[pdftex,pdftitle={anexoc},colorlinks=true,breaklinks=true,linkcolor=blue]{hyperref}
\font\cmbdoce=cmb10 at 12pt

\font\cmbcatorce=cmb10 at 14pt
\font\cmronce=cmr10 at 11pt
\font\cmita=cmti10 at 12pt
\newcommand{\supp}{\textsf{supp}}
\begin{document}
\vspace*{14pt}
\begin{center}
{  \cmbcatorce  MULTIRESOLUTION SCHEMES AND ITS APPLICATION TO SEDIMENTATION MODELS}

\vspace{14pt}

{\cmbdoce Ricardo Ruiz-Baier}
\vspace{6pt}

{\cmronce Departamento de Ingenier\'{\i}a Matem\'{a}tica,

Universidad de Concepci\'{o}n,

Casilla 160-C, Concepci\'{o}n, Chile.

e-mail: rruiz@ing-mat.udec.cl

}

\vspace{12pt}
\end{center}

\vspace{12pt}
\noindent {\cmbdoce Palabras clave:} M\'etodos de Multiresoluci\'on, Wavelets, Leyes de Conservaci\'on, Procesos de Sedimentaci\'on, Problemas de flujo vehicular

\vspace{12pt}

\noindent {\cmbdoce Resumen.}{ \cmita
Se presenta un m\'etodo num\'erico para obtener soluciones aproximadas de problemas provenientes de
la sedimentaci\'on de suspensiones floculadas. Estos procesos se usan ampliamente en la industria minera, por ejemplo para recuperar el agua de
las suspensiones que salen de los procesos de flotaci\'on \cite{Thick}.

La idea principal es aplicar los m\'etodos de multiresoluci\'on a los esquemas desarrollados por B\"urger \emph{et al.}
\cite{BCS,BEKL,BK} y observar que el m\'etodo de multiresoluci\'on a describir es vital para reducir el costo computacional sin afectar la calidad de la soluci\'on.

\vspace{12pt}

}

\begin{multicols}{2}

\vspace{12pt}
\noindent{\cmbdoce 1 INTRODUCCION}
\vspace{6pt}

Se introduce el problema f{\'i}sico y su modelaci\'on mediante una ley de conservaci\'on fuertemente degenerada
con flujo no lineal. El efecto de la compresibilidad del sedimento puede ser descrito por un t\'ermino difusivo fuertemente
degenerado, mientras el flujo unidimensional contribuye una discontinuidad de flujo a la ecuaci\'on parcial diferencial. Se presenta un esquema de segundo orden para resolver este tipo de problemas y finalmente se desarrollan ejemplos num\'ericos para comparar con resultados publicados (\cite{BEKL,BK}).

En el caso unidimensional, la teor{\'i}a de la sedimentaci\'on produce ecuaciones de equilibrio de masa y momentum lineal que pueden simplificarse hasta obtener una ecuaci\'on de la forma
\begin{equation}\label{degen1}
\partial_tu+\partial_xf(u)=\partial_{xx}^2A(u),
\end{equation}
con $(x,t)\in ]0,1[\times[0,T[$ y el \emph{coeficiente de difusi\'on integrado} dado por
\begin{equation}
A(u):=\int_0^ua(s)ds,\quad a(u)\geqslant 0.
\end{equation}
Se permite que el \emph{coeficiente de difusi\'on} $a(u)$ sea cero sobre intervalos de $u$.

Las soluciones de (\ref{degen1}) desarrollan discontinuidades debido a la no linealidad de la funci\'on de densidad de flujo $f(u)$ y a la degeneraci\'on del coeficiente de difusi\'on. Esto lleva a considerar soluciones entr\'opicas para tener un problema bien puesto. A\'un m\'as, cuando (\ref{degen1}) es puramente hiperb\'olica, los valores de la soluci\'on se propagan sobre rectas caracter{\'i}sticas que podr{\'i}an intersectar las fronteras del dominio espacio-tiempo desde el interior, y esto requiere tratar a las condiciones de Dirichlet como condiciones entr\'opicas.

$a(u)$ tiene un \emph{comportamiento degenerado}, es decir, $a(u)=0$ para $u\leqslant u_c$ y $a(u)$ salta en $u_c$ a un valor positivo, donde $u_c$ es una constante llamada \emph{concentraci\'on cr{\'i}tica}. Se enfatiza el hecho de que el coeficiente de difusi\'on $a(u)$ es degenerado, lo que hace evidente la naturaleza hiperb\'olica-parab\'olica de (\ref{degen1}).

Considerar el PVIF siguiente
{\tiny
\begin{eqnarray}
\partial_tu+\partial_x(q(t)u+f(u))&=&\partial_{xx}^2A(u),\quad (x,t)\in ]0,H[\times[0,T[,\label{A1}\\
u(x,0)&=&u_0(x),\quad x\in [0,H],\label{A2}\\
u(H,t)&=&0,\quad t\in]0,T]\label{A3}\\
f(u(0,t))-\partial_xA(u(0,t))&=&0,\quad t\in]0,T],\label{A4}
\end{eqnarray}}
conocido como el \emph{Problema A}.
Y el \emph{Problema B}
{\tiny\begin{eqnarray}
\partial_tu+\partial_x(q(t)u+f(u))&=&\partial_{xx}^2A(u),\quad (x,t)\in ]0,H[\times[0,T[,\label{B1}\\
u(x,0)&=&u_0(x),\quad x\in [0,H],\label{B2}\\
q(t)u(H,t)-\partial_xA(u(H,t))&=&\Psi(t),\quad t\in]0,T]\label{B3}\\
f(u(0,t))-\partial_xA(u(0,t))&=&0,\quad t\in]0,T].\label{B4}
\end{eqnarray}}
Para ambos problemas, $f$ se supone continua y diferenciable a trozos, $f\leqslant 0$, $\supp(f)\subset [0,u_{\max}]$, $\|f'\|_\infty\leqslant \infty$, $a(u)\geqslant 0$, $\supp(a)\subset \supp(f)$, $a(u)=0$ para $u\leqslant u_c$, $0<u_c<u_{\max}$, $q(t)\leqslant 0,\ \forall t\in [0,T]$, $TV(q)<\infty,\ TV(q')<\infty$. En \cite{BEKL} se prueba la existencia y unicidad de soluci\'on entr\'opica para cada uno de estos problemas.

La propiedad de mayor inter\'es, es que generalmente se supone el siguiente comportamiento para $\sigma_e(u)$:
{\tiny\begin{equation}
\sigma_e(u)\left\{\begin{array}{ll}
=\textrm{cte.},& \textrm{ si } u\leqslant u_c,\\
>0,& \textrm{ si } u>u_c,\end{array}\right.\textrm{ y }\sigma_e'(u):=\frac{d\sigma_e}{du}\left\{\begin{array}{ll}
=0,& \textrm{ si } u\leqslant u_c,\\
>0,& \textrm{ si } u>u_c.\end{array}\right.
\end{equation}}

\vspace{12pt}
\noindent{\cmbdoce 2 ESQUEMAS DE SEGUNDO ORDEN}
\vspace{6pt}

Los t\'erminos advectivo y difusivo son aproximados de diferente forma, con el fin de obtener una discretizaci\'on que mantenga la conservatividad. Para la parte advectiva puede utilizarse el esquema de Roe cl\'asico con una interpolaci\'on ENO de segundo orden, o bien un esquema de Engquist-Osher \cite{EO}. Para la parte difusiva, se necesita un esquema centrado de segundo orden que mantenga la conservatividad. El esquema interior resultante para la ecuaci\'on (\ref{degen1}) es
{\tiny\begin{eqnarray*}
&&\frac{u_j^{n+1}-u_j^n}{\Delta t}+q(n\Delta t)\frac{u^-_{j+1}-u^+_{j-1}}{\Delta x}+\frac{F_{j+\frac{1}{2}}-F_{j-\frac{1}{2}}}{\Delta x}=\\
&&\frac{A(u_{j-1}^n)-2A(u_j^n)+A(u_{j+1}^n)}{(\Delta x)^2}
\end{eqnarray*}}
En este caso, se utilizar\'a el \emph{$\theta-$limitador}
{\tiny $$ s_j^n=MM\left(\theta\frac{u_j^n-u_{j-1}^n}{\Delta x},\, \frac{u_{j+1}^n-u_{j-1}^n}{2\Delta x},\, \theta\frac{u_{j+1}^n-u_j^n}{\Delta x}\right),\quad \theta\in[0,2],$$}
con
{\scriptsize \begin{equation}
MM(a,b,c):=\left\{\begin{array}{ll}
\min(a,b,c),&\textrm{ si } a,b,c>0,\\
\max(a,b,c),&\textrm{ si } a,b,c<0,\\
0,&\textrm{e.o.c.}\end{array}\right.
\end{equation}}
Este esquema es estable bajo la condici\'on $CFL$
{\small\begin{equation}
\max_u|f'(u)|\frac{\Delta t}{\Delta x}+2\max_u|a(u)|\frac{\Delta t}{(\Delta x)^2}\leqslant 1.
\end{equation}}

\vspace{12pt}
\noindent{\cmbdoce 3 ANALISIS DE MR}
\vspace{6pt}

Se presentan los conceptos y definiciones b\'asicas introducidas por Harten \cite{Harten1} para el an\'alisis de multiresoluci\'on.
Considerar el conjunto de mallas anidadas di\'adicas $G^k,\ k=0,\ldots,L$:
{\scriptsize\begin{equation}
G^k=\{x_j^k\}_{j=0}^{N_k},\ x_j^k=-1+j\cdot h_k,\ h_k=2^{N_k+1}h_0,\ N_k=\frac{N_0}{2^k},
\end{equation}}
\begin{center}
\includegraphics[height=3cm]{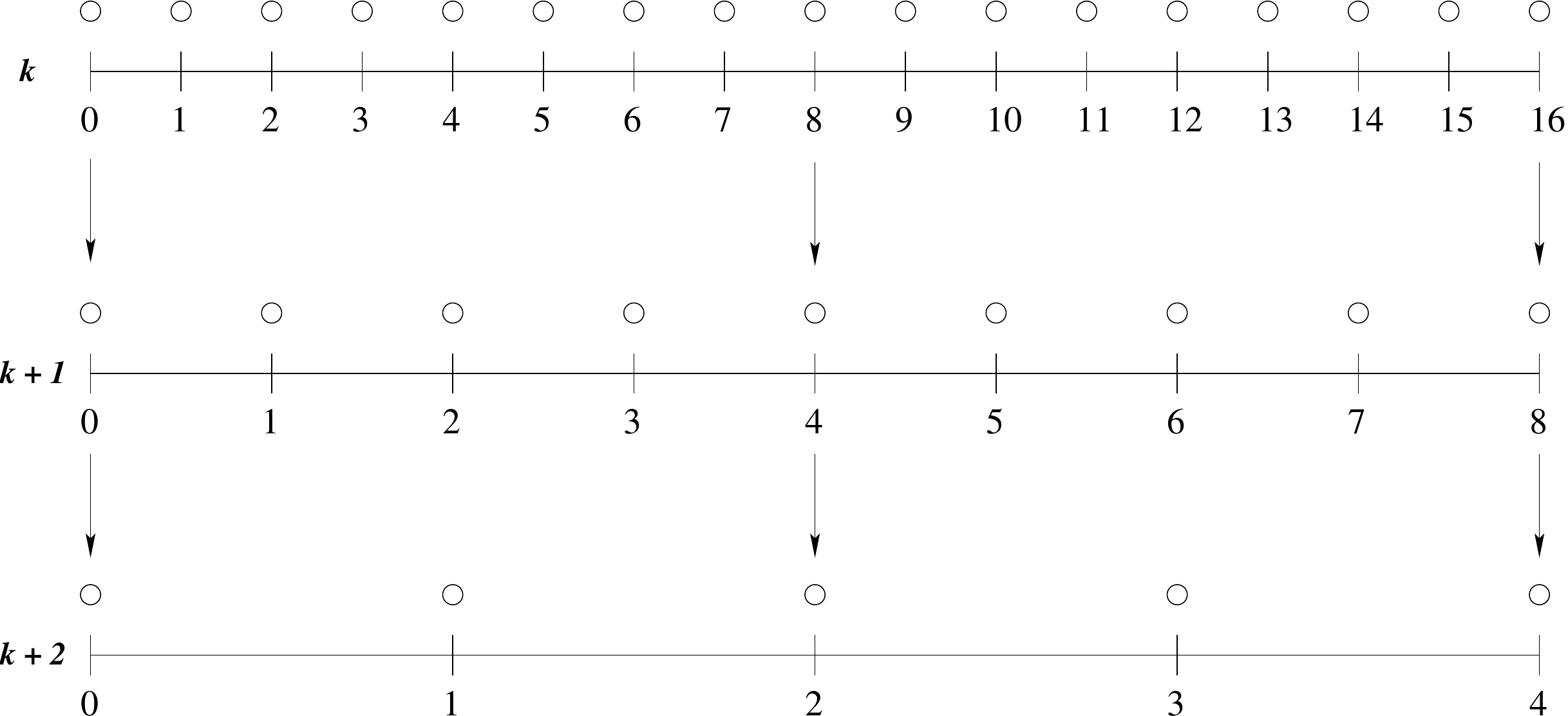}

{\scriptsize Diferentes escalas de valores puntuales}
\end{center}

\vspace{12pt}
\noindent{\cmbdoce 4 ALGORITMO DE MR}
\vspace{6pt}

\begin{enumerate}
\item \emph{Inicializaci\'on de par\'ametros y variables}:
\begin{itemize}
\item Longitud del dominio $H$, concentraci\'on cr{\'i}tica $u_c$, orden de la interpolaci\'on de multiresoluci\'on $r$, niveles de multiresoluci\'on $L$, n\'umero de puntos y paso en la malla fina $N_0$ y $h_0$, y en cada nivel $N_k$ y $h_k$, tolerancia prescrita $\varepsilon$ y estrategia de truncamiento $\varepsilon_k$, constantes de Lipschitz para $a(u)$ y $f'(u)$,
\item condici\'on $CFL$:
$$\max_u|f'(u)|\frac{\Delta t}{h_0}+2\max_u|a(u)|\frac{\Delta t}{h_0^2}\leqslant 1.$$
\item paso temporal $\Delta t$,
$$\Delta t=\frac{CFL\cdot h_0}{\max_u|f'(u)|+2\max_u|a(u)|/h_0}.$$
\item estructura de datos. (SPARSE, GRADED TREE).
\end{itemize}
\item \emph{Aplicaci\'on de la codificaci\'on a la condici\'on inicial}:
Coeficientes de ondelette significativos y posiciones correspondientes. Se incluyen los \emph{safety points}.
\item \emph{Evoluci\'on temporal}: Se utiliza un m\'etodo
Runge-Kutta de segundo orden.
\end{enumerate}

\vspace{12pt}
\noindent{\cmbdoce 5 EJEMPLOS NUMERICOS }
\vspace{6pt}

Se calculan soluciones de los problemas A y B utilizando los
esquemas num\'ericos descritos en la secci\'on anterior, con
una discretizaci\'on para el flujo de tipo Enqguist-Osher. Se
reproducen algunos resultados num\'ericos.

\vspace{12pt}
\noindent{\cmbdoce 5.1 SEDIMENTACION BATCH IDEAL }
\vspace{2pt}

Proceso de sedimentaci\'on batch de suspensi\'on ideal en una columna de asentamiento:
\begin{center}
\includegraphics[width=7.5cm,height=5cm]{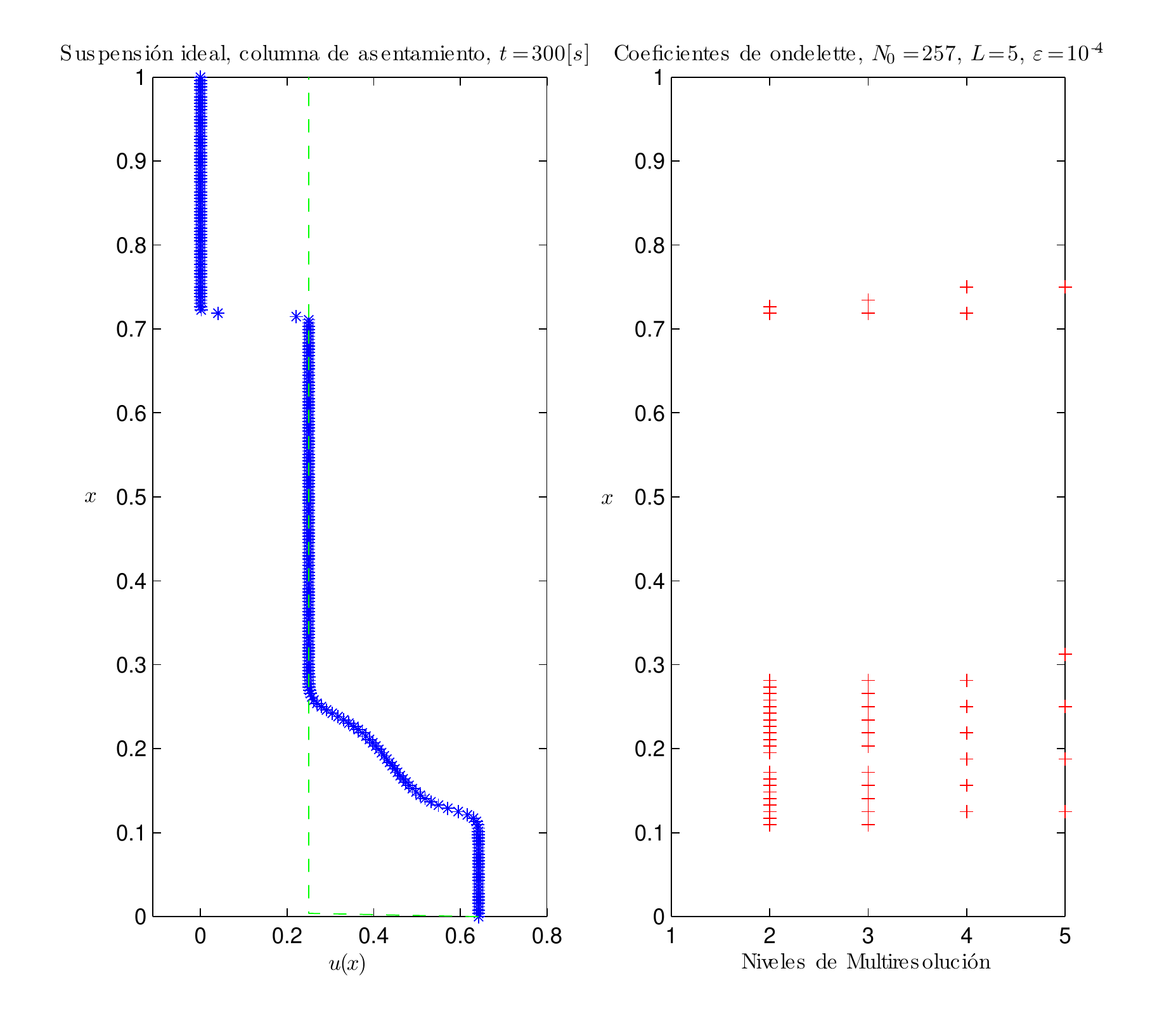}

{\tiny Izquierda: Condici\'on inicial \emph{(rayas)} y perfil de concentraci\'on a $t=300[s]$ para el problema de sedimentaci\'on  batch de suspensi\'on ideal \emph{(Asteriscos)}. Derecha: Coeficientes de ondelette significativos correspondientes.}
\end{center}

{\tiny
\begin{center}
\begin{tabular}{lcccc}
\hline
$t\, [s]$ & $V$    & $\mu$  &            $e_1$   &    $e_\infty$     \\
\hline
60  &4.3457 & 7.8456 & 2.64$\times10^{-5}$&9.03$\times10^{-6}$\\
300 &5.6212 & 5.8456 & 1.70$\times10^{-5}$&1.12$\times10^{-5}$\\
1800&5.9443 & 14.9168 & 7.28$\times10^{-5}$&4.35$\times10^{-5}$\\
3600&6.1385 & 29.8479& 8.89$\times10^{-5}$&6.50$\times10^{-5}$\\
\hline
\end{tabular}
\end{center}}
{\tiny Sedimentaci\'on de suspension ideal. $\varepsilon=1.0\times10^{-4}$, $N_0=257$ y $L=5$.}

\vspace{6pt}
\noindent{\cmbdoce 5.2 SEDIMENTACION BATCH CON COMPRESION }
\vspace{2pt}

Como funci\'on de densidad de flujo, se utiliza una funci\'on Kynch batch Richardson-Zaki con par\'ametros correspondientes a suspensi\'on de cobre \cite{BEKL}.
\begin{equation}\label{zaki}
f(u)=-6.05\times10^{-4}u(1-u)^{12.59}\, [m/s].
\end{equation}
Se utilizar\'a la funci\'on $\sigma_e'(u)$ dada por (\cite{BK,Thick})
\begin{equation}
\sigma_e'(u)=\frac{d}{d\,u}\left(100(u/u_c)^8-1\right)[Pa], \textrm{ si } u> u_c.
\end{equation}
Luego
\begin{equation}
\sigma_e'(u)=\left\{\begin{array}{ll}
0,&\textrm{ si } u\leqslant u_c=0.23,\\
\frac{800}{u_c}\left(\frac{u}{u_c}\right)^7\, [Pa],& \textrm{ si } u> u_c.
\end{array}\right.
\end{equation}

\begin{center}
{\tiny
\begin{tabular}{lcccc}
\hline
$t\, [s]$ & $V$    & $\mu$  &            $e_1$   &       $e_\infty$     \\
\hline
60       &6.5737 & 17.8796& 1.29$\times10^{-4}$&5.33$\times10^{-5}$\\
1800 (*) &5.7349 & 9.4132 & 1.99$\times10^{-4}$&7.42$\times10^{-5}$\\
3600 (*) &6.1982 & 9.1246 & 2.77$\times10^{-4}$&9.61$\times10^{-5}$\\
7200 (*) &6.2110 & 9.1246 & 3.21$\times10^{-4}$&2.41$\times10^{-4}$\\
14400(*) &7.9244 & 9.4132 & 8.92$\times10^{-4}$&6.18$\times10^{-4}$\\
\hline
\end{tabular}}

{\tiny Suspensiones floculadas, primer ejemplo.
Multiresoluci\'on utilizando $\varepsilon=10^{-3}$,
$N_0=129$ y $L=5$.}
\end{center}
\begin{center}
\includegraphics[width=7.5cm,height=5cm]{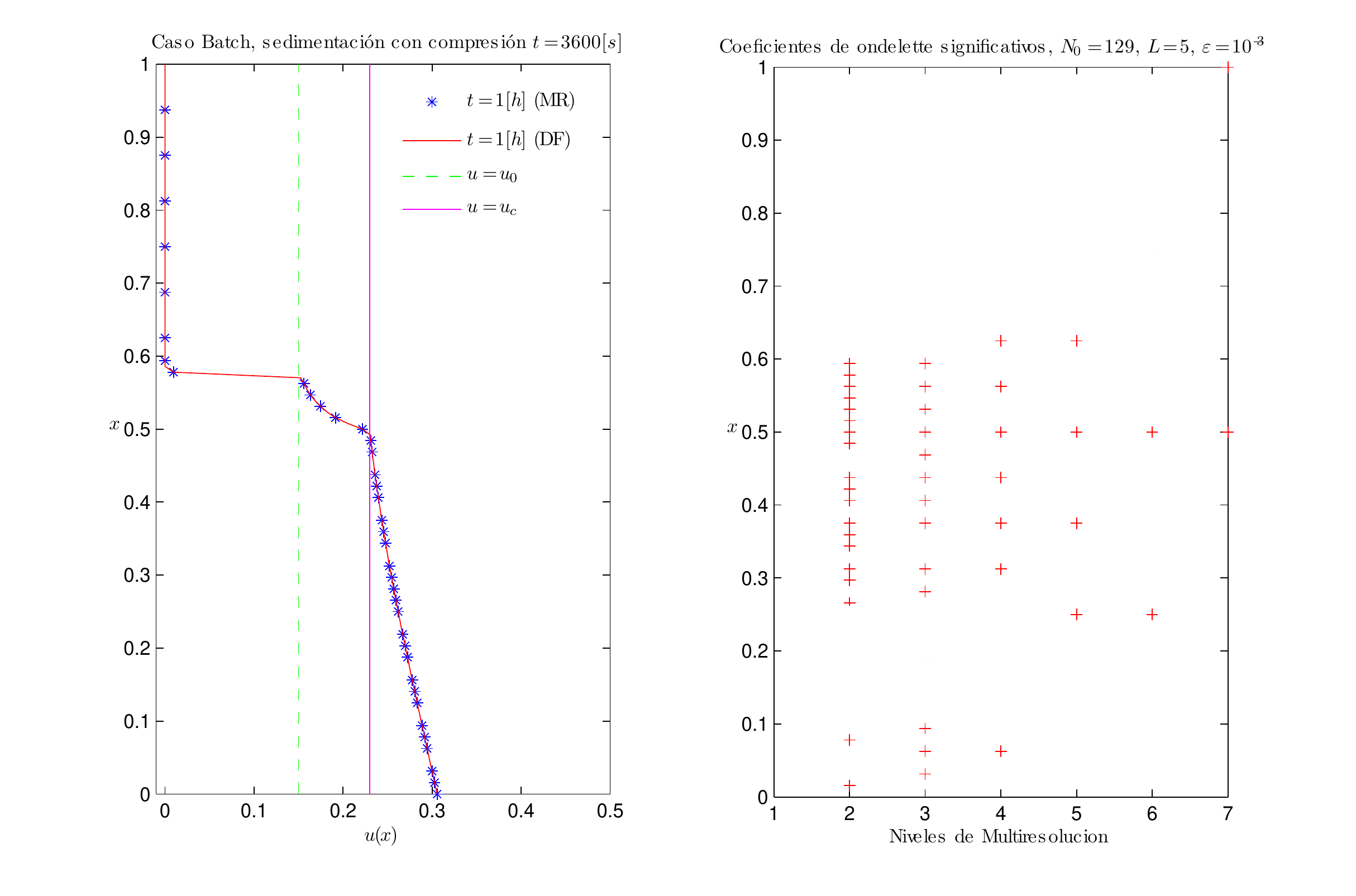}
{\tiny Izquierda: Condici\'on inicial \emph{(rayas)} y perfil de
concentraci\'on a $t=3600[s]$ para el problema de
sedimentaci\'on-consolidaci\'on  \emph{(asteriscos)}. Derecha:
Coeficientes de ondelette significativos correspondientes.
$\varepsilon=10^{-3}$, $N_0=129$ y $L=5$.}
\end{center}
\begin{center}
\includegraphics[width=7.5cm,height=3.8cm]{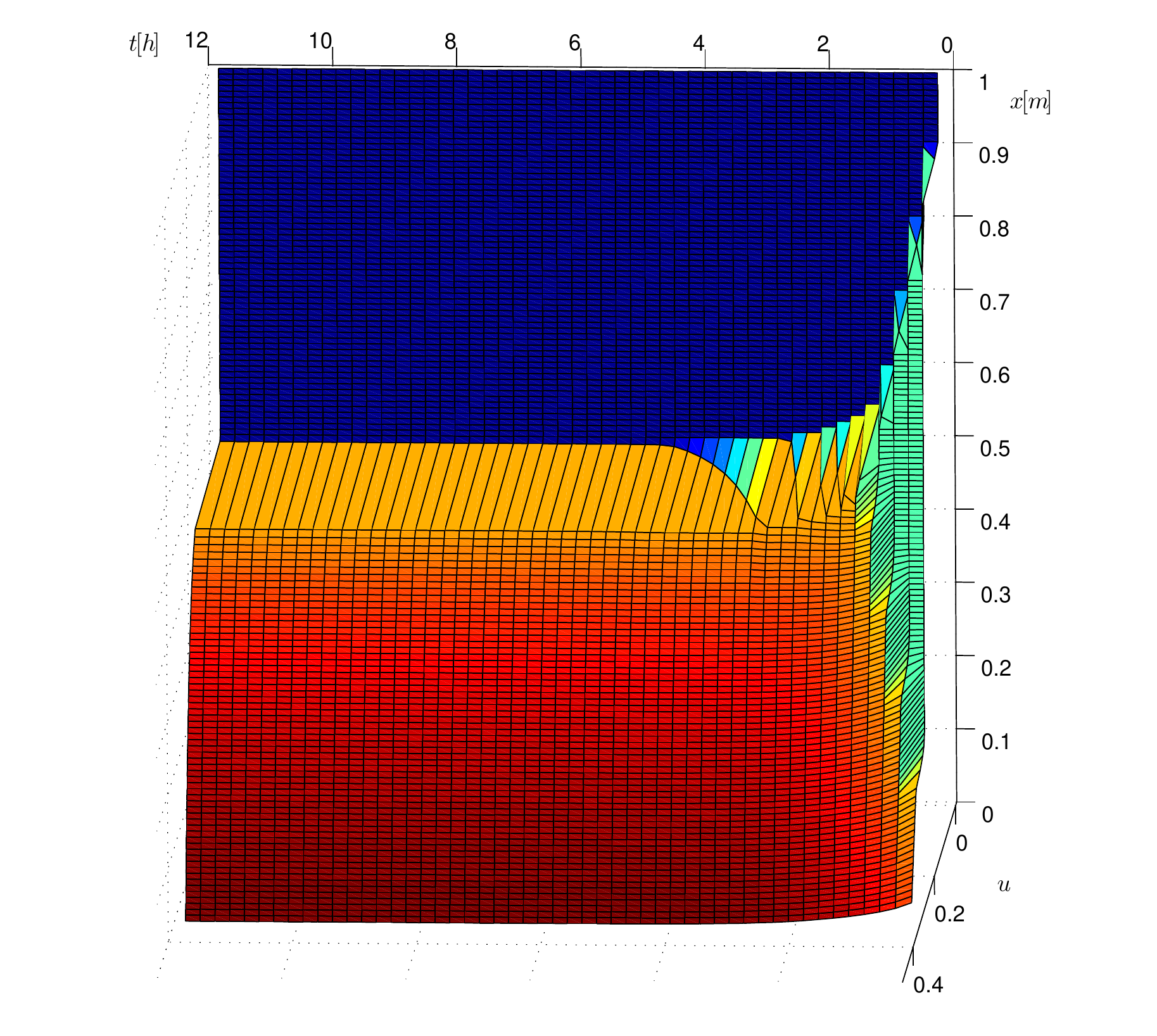}
{\tiny Perfiles de concentraci\'on hasta $t=12[h]$, problema de sedimentaci\'on-consolidaci\'on \emph{Batch.}  $\varepsilon=10^{-3}$, $N_0=129$ y $L=5$.}
\end{center}
\vspace{2pt}
\noindent{\cmbdoce 5.3 SEDIMENTACION CONTINUA }
\vspace{2pt}
\begin{center}
\includegraphics[width=5cm,height=3.5cm]{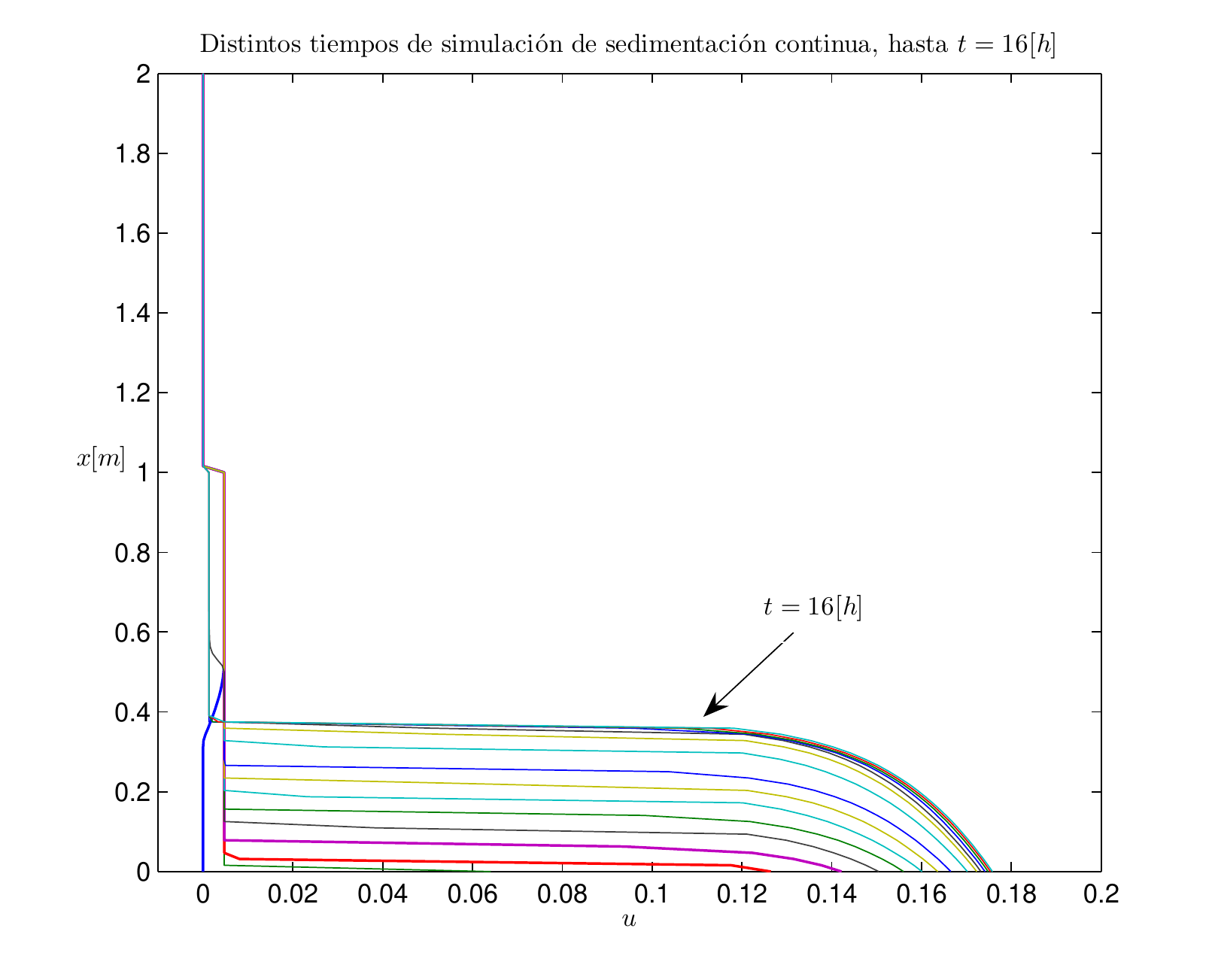}

{\tiny Perfiles de concentraci\'on hasta
$t=16[h]$ para el problema de sedimentaci\'on continua. $\varepsilon=5\times10^{-4}$, $N_0=513$ y $L=5$.}
\end{center}
\vspace{2pt}
\noindent{\cmbdoce 5.4 REACCION-DIFUSION 1D }
\vspace{1pt}
\begin{center}
\includegraphics[width=6cm,height=3cm]{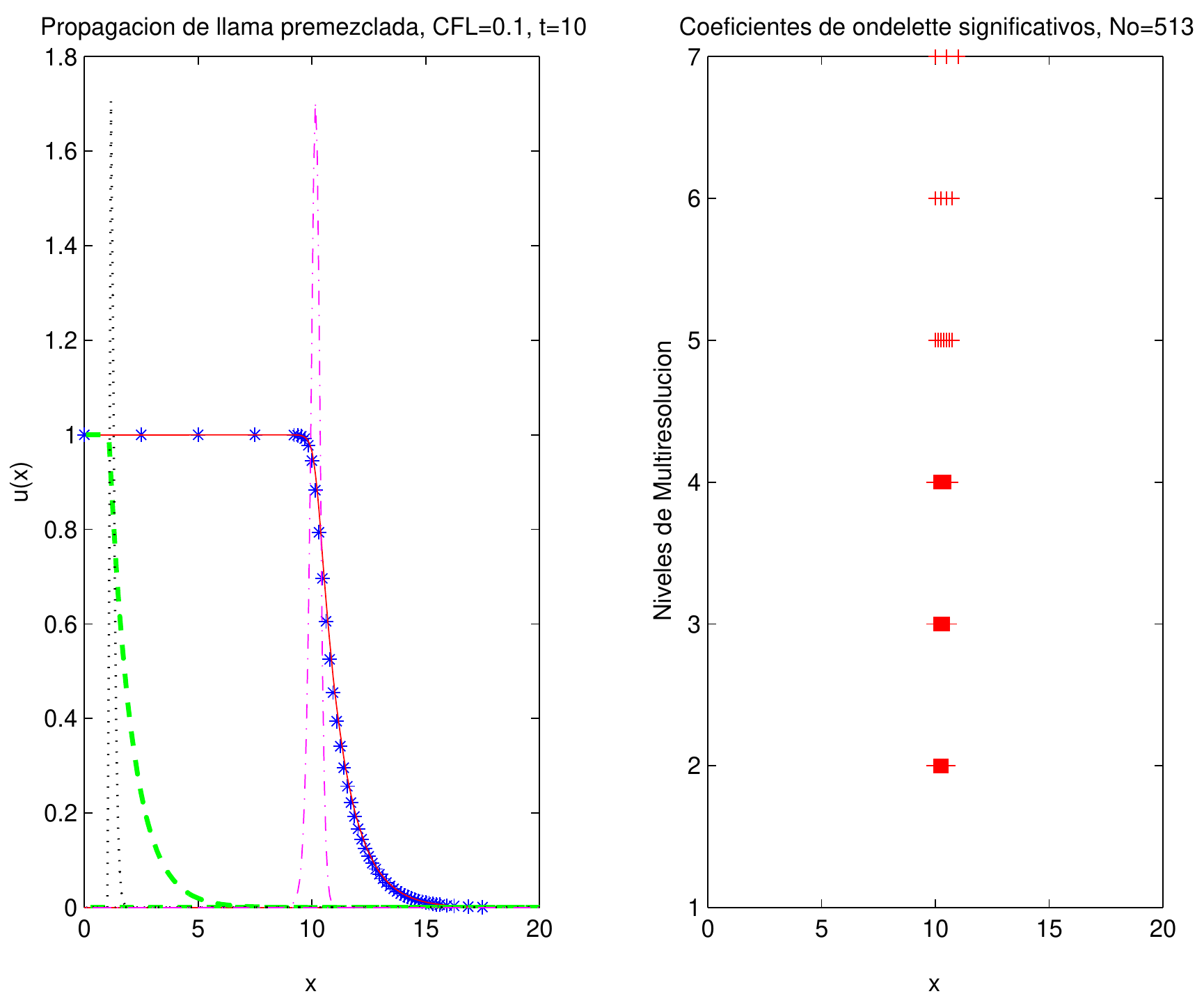}

{\tiny Condici\'on inicial \emph{(rayas)} y $S(u)$ inicial \emph{(puntos)}, soluci\'on num\'erica sin multiresoluci\'on \emph{(linea)}, soluci\'on num\'erica con multiresoluci\'on \emph{(asteriscos)} y $S(u)$ \emph{(puntos-rayas)}, en el tiempo $t=10$ para la ec. de reacci\'on-difusi\'on, con $\alpha=0.8$, $\beta=10$, $L=7$, $N_0=513$ y $\varepsilon=10^{-3}$. Detalles significativos, $t=0.5$.}
\end{center}
\vspace{1.2pt}
\noindent{\cmbdoce 5.5 TRAFICO VEHICULAR }
\vspace{-2pt}
\begin{center}
\includegraphics[width=0.47\textwidth]{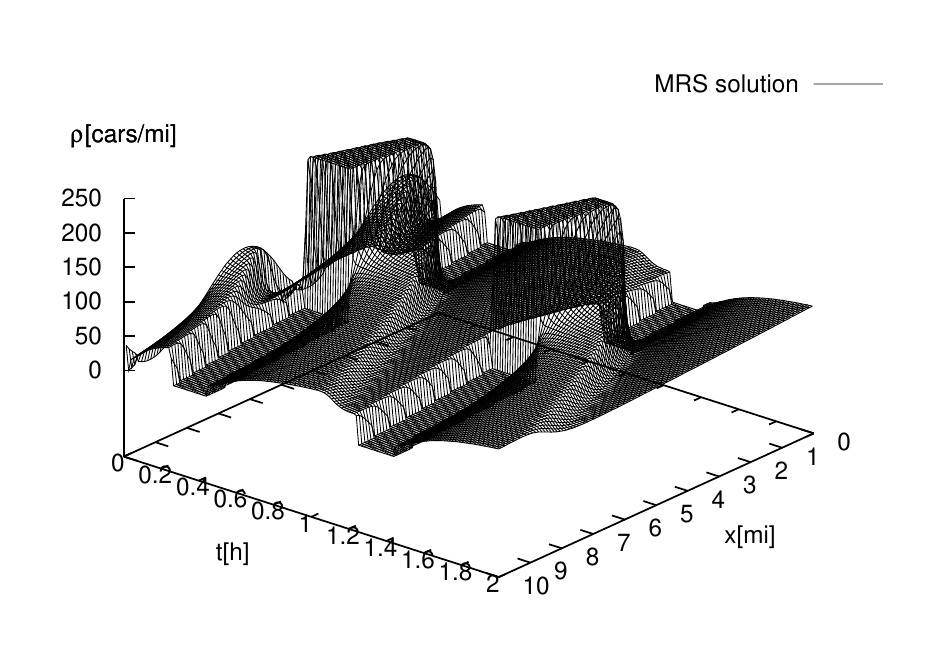}
{\tiny Soluci\'on tridimensional del problema de flujo de tr\'afico en una rotonda.}

\includegraphics[width=0.47\textwidth]{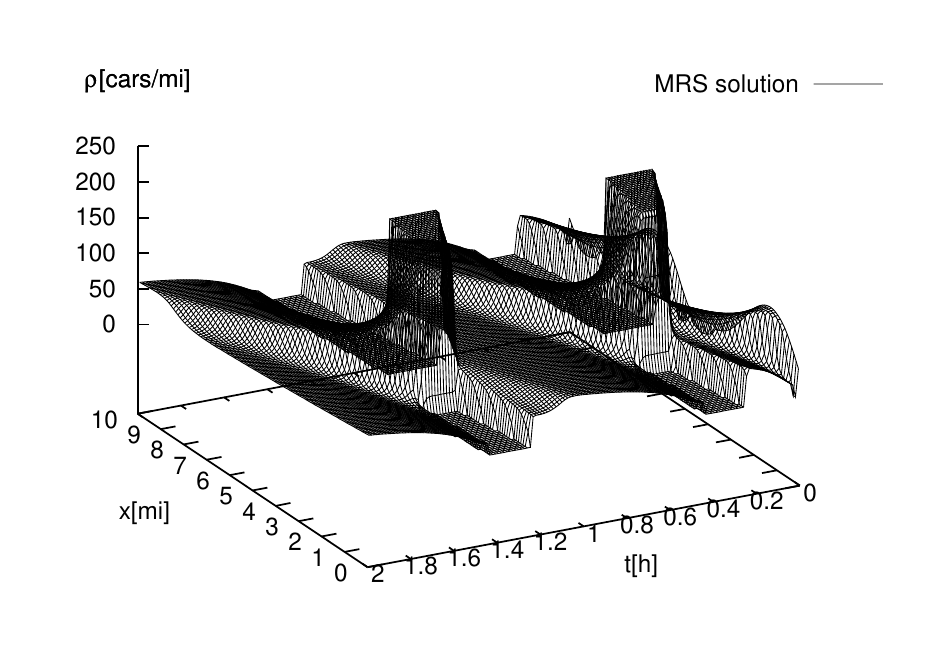}
{\tiny Soluci\'on tridimensional del problema de flujo de tr\'afico en una rotonda.}
\end{center}

\end{multicols}

\end{document}